\documentclass[a4paper, 11pt]{amsart} 

\usepackage{amsmath,amsthm,amssymb,mathtools,url,breqn}

\usepackage{ulem}
\usepackage{color}
\usepackage{amssymb}
\usepackage{amstext}
\usepackage{amsmath}
\usepackage{amscd}
\usepackage{amsthm}
\usepackage{amsfonts}
\usepackage{enumerate}
\usepackage{graphicx}
\usepackage{color}
\usepackage{here} % 図や表を強制的に出力させる
\usepackage{bm} % 太字ベクトル
\usepackage{mathrsfs}
\usepackage{mathtools}
\usepackage{latexsym}
\usepackage{booktabs}
\usepackage{fancyhdr}
\usepackage{subcaption}
\usepackage{tikz}
\usepackage{tikz-cd}
\usetikzlibrary{arrows}
\usetikzlibrary{decorations.markings}
\usetikzlibrary{patterns,decorations.pathreplacing}

\newcommand{\PF}{\textrm{PF}}

\newcommand{\Ape}{\mathrm{Ap}}

\newcommand{\move}{\textrm{M}}

\newcommand{\tr}{\textrm{tr}}
\newcommand{\bfa}{\mathbf{a}}
\newcommand{\mm}{\mathbf{m}}

\newcommand{\bfv}{\mathbf{v}}
\newcommand{\bfu}{\mathbf{u}}

\newcommand{\bfx}{\mathbf{x}}

\newcommand{\ZZ}{\mathbb{Z}}
\newcommand{\NN}{\mathbb{N}}

\newcommand{\PP}{\mathbb{P}}

\newcommand{\AAA}{\mathbb{A}}

\newcommand{\kk}{\Bbbk}

\newtheorem{mainthm}{Theorem}

\newtheorem{mainlem}{Lemma}

\theoremstyle{definition}
\newtheorem{dfn}{Definition}[section]
\newtheorem{rem}[dfn]{Remark}

\newtheorem{ex}[dfn]{Example}

\theoremstyle{plain}

\newtheorem{thm}[dfn]{Theorem}
\newtheorem{lem}[dfn]{Lemma}
\newtheorem{prop}[dfn]{Proposition}
\newtheorem{cor}[dfn]{Corollary}

     \title{Nearly Gorenstein projective monomial curves of small codimension}

     \author{Sora Miyashita}
     \address{Department of Pure And Applied Mathematics, Graduate School Of Information Science And Technology, Osaka University, Suita, Osaka 565-0871, Japan}
     \email{u804642k@ecs.osaka-u.ac.jp}
     \date{October, 4, 2023}


\begin{document}
     \begin{abstract}
     In this paper, we characterize nearly Gorenstein projective monomial curves of codimension 2 and 3.
     \end{abstract}
     \maketitle

\section{Introduction}
Let $\kk$ be a field,
and let us denote the set of nonnegative integers and the set of integers by $\mathbb{N}$ and $\mathbb{Z}$, respectively.
Cohen-Macaulay rings and Gorenstein rings are important properties and play a crucial role in the theory of commutative algebras.
Many kinds of rings are defined for studying of a new class of local or graded rings which are Cohen-Macaulay but not Gorenstein.
For example, there are \textit{nearly Gorenstein} rings, \textit{almost Gorenstein} rings and
\textit{level} rings, and so on.
Let $R$ be a Cohen-Macaulay $\mathbb{N}$-graded $\kk$-algebra with canonical module $\omega_R$.
According to \cite{HHS}, $R$ is called nearly Gorenstein if the trace ideal {of the canonical} module $\tr(\omega_R)$ contains the maximal graded ideal $\mm$ of $R$.
Here, let $\tr(\omega_R)$ be the ideal
generated by the image of $\omega_R$ through all homomorphism of $R$-modules into $R$.
In particular, $R$ is non-Gorenstein and nearly Gorenstein
if and only if $\tr(\omega_R) = \mathbf{m}$ ({see \cite[Definition 2.2]{HHS}}).
On the other hand, according to \cite[Chapter III, Proposition 3.2]{Sta},
$R$ is called level if the degrees of the minimal
generators of $\omega_R$ are the same.

There are characterizations of nearly Gorenstein rings on some concretely graded rings.
For example, it is known on Hibi rings, Stanley-Reisner rings with Krull dimension 2 and numerical semigroup rings  with small embedding dimension, and so on (see \cite{nume,HHS,M}).
In this paper, we characterize nearly Gorenstein projective monomial curves of codimension 2 and 3.
Let $\PP_{\kk}$ be a projective space over $\kk$, and let $C$ be a monomial curve in $\PP_{\kk}^n$.
We consider its projective coordinate ring $\AAA(C)$.
In fact,
this ring is {isomorphic} to
the affine semigroup ring $\kk[S_{\bfa}]$ where
$$S_\bfa =  \langle(0,a_n),(a_1,a_n-a_1),(a_2,a_n-a_2),\cdots,
(a_{n-1},a_n-a_{n-1}),(a_n,0) \rangle.$$
We can assume $0<a_1<a_2<\cdots<a_n$ are integers
with $\gcd(a_1,a_2,\cdots,a_n)=1$.
We call this semigroup $S_\bfa$
{the projective} monomial curve defined by $\bfa$.
We call $S_{\bfa}$ is Cohen-Macaulay, Gorenstein,
nearly Gorenstein and level if its
affine semigroup ring $\kk[S_\bfa]$ is so, respectively.
If $S_{\bfa}$ is Cohen-Macaulay, we denote by $r(S_\bfa)$ the Cohen-Macaulay type
of $\kk[S_{\bfa}]$.
In Proposition \ref{4.2} and \ref{4.3},
we prove the following result. This is our main result.

\begin{mainthm}\label{thm:mainA}
Let $\bfa=a_1,\cdots,a_n$ be a sequence of positive integers with
$a_1<a_2<\cdots<a_n$ and $\gcd(a_1,\cdots,a_n)=1$.
Let $S_\bfa$ be the projective monomial curve defined by $\bfa$.
Assume that $S_{\bfa}$ is Cohen-Macaulay.
\begin{itemize}
\item If $n=3$,
  then the following conditions are equivalent:
\begin{itemize}
\item[(1)] $S_{\bfa}$ is non-Gorenstein and nearly Gorenstein;
\item[(2)] $\bfa=k,k+1,2k+1$
  for some $k \geqq 1$.
\end{itemize}
If this is the case, $S_{\bfa}$ is level with $r(S_{\bfa})=2$.
  \item If $n=4$,
then the following conditions are equivalent:
\begin{itemize}
\item[(1)] $S_{\bfa}$ is non-Gorenstein and nearly Gorenstein;
\item[(2)] $\bfa=1,2,3,4$ or $S_\bfa \cong S_{2k-1,2k+1,4k,6k+1}$ for some $k \geqq 1$.
\end{itemize}
If this is the case, $S_{\bfa}$ is level with $r(S_{\bfa})=3$.
\end{itemize}
\end{mainthm}

In particular, we get the following corollary from Theorem A.

\begin{cor}
Every nearly Gorenstein projective monomial curve of codimension at most 3 is level.
\end{cor}

\begin{rem}
{If $n=2$, then $\kk[S_{\bfa}]$ is always Gorenstein because it is hypersurface.}
\end{rem}

\begin{rem}
It is known that there exists a nearly Gorenstein but not level projective monomial curve of codimension 4.
Indeed, $S_{4,9,12,13,21}$ is such a example
({see \cite[Theorem 3.11]{M}}).
\end{rem}

The structure of this paper is as follows.
In Section 2, we prepare some facts and lemmas to
show Theorem A.
In Section 3, we discuss nearly Gorensteinness of some concretely projective monomial curves
to prove Theorem A.
In Section 4, we prove Theorem A.

\subsection*{Acknowledgement}
I am grateful to Professor
Akihiro Higashitani for his very
helpful comments and instructive discussions.
I also thank Tomoya Hikida
for his help in implementing
the program to observe projective monomial curves.
\section{Preliminaries}

\subsection{Homogeneous affine semigroup rings}
An \textit{affine semigroup} $S$ is a finitely generated sub-semigroup of $\mathbb{N}^d.$
An affine semigroup $S$ is \textit{homogeneous} if all its minimal generators lie on an affine hyperplane not including origin.
This is equivalent to the condition {that the affine} semigroup ring $\kk[S]$ is standard graded by assigning degree $1$ to all the monomials corresponding to the minimal generators of $S$.
We denote by $\ZZ S$ the group generated by $S$.
%In that case, we also say that $\kk[S]$ is %\textit{homogeneous}.
Let $S$ be a Cohen-Macaulay homogeneous affine semigroup,
and let $G_S=\{ \bfa_1,\cdots,\bfa_s \} \subseteq \mathbb{N}^n$ be the minimal generators of $S$.
Fix the affine semigroup ring $R=\kk[S]$.
Since $S$ is homogeneous, we can regard
$R=\kk[S]=\kk[\bfx^{\bfa_1},\cdots,\bfx^{\bfa_s}]$ as standard graded by assigning
$\deg \bfx^{\bfa_i}=1$
for any $1 \leq i \leq s$.
In this case, the canonical module $\omega_{R}$ is isomorphic to an ideal $I_{R}$ of $R$
as an {$\mathbb{N}^n$}-graded module up to degree shift.
Then we can assume the system of minimal generators of $I_R$
is $\{ \bfx^{\bfv_1},\cdots,\bfx^{\bfv_r} \}$,
and $V(S)=\{\bfv_1,\cdots,\bfv_r\} \subseteq S$ is a minimal {generating system of the} canonical ideal of $S$.
We put
$V_{\min}(S)
=
\{\bfv \in V(S)
; \deg \bfx^\bfv \leqq
\deg \bfx^{\bfv_i} \;\text{for all}
\; 1 \leqq i \leqq r\}$
and
$S-V(S)
:= \{ \bfa \in \mathbb{Z}S \;;\; \bfa + \bfv \in S \; \text{for all}\; \bfv \in V(S)\}$.
Thus the following holds.
\begin{prop}[{see \cite[Proposition 3.4]{M}}]\label{2.1}
  Let $S$ be a Cohen-Macaulay
  homogeneous affine semigroup.
  The following are equivalent:
  \begin{itemize}
\item[(1)] $R=\kk[S]$ is nearly Gorenstein;
\item[(2)] For any $\bfa_i \in G_S$,
  there exist $\bfv \in V_{\min}(S)$
  and $\bfu \in S-V(S)$
  such that
  $\bfa_i = \bfu + \bfv$.
  \end{itemize}
\end{prop}

For any set $X$, we denote by $|X|$ the cardinality of $X$.
Now we discuss $|V_{\min}(S)|$ to show Theorem A.

\begin{thm}[{see \cite[Theorem 3.6]{M}}]
  If $R$ is non-Gorenstein and nearly Gorenstein,
  then $|V_{\min}(S)| \geqq 2$.
  \end{thm}

%\begin{dfn}[{see \cite{teter}}]
%  \textit{Teter number} of $R$ to be the smallest
%  number $s$ for which there exist $R$-module homomorphisms
%  $\phi_i : \omega_R \rightarrow R$ such that
%  $\tr(\omega_R) = \sum_{i=1}^s \phi_i(\omega_R).$
%  Non-Gorenstein ring $R$ is \textit{Teter type} if
%  its Teter number is 1.
%\end{dfn}

%\begin{thm}[{see \cite[Theorem 1.1]{teter}}]
%Let $R$ be a non-Gorenstein ring.
%If $R$ is a domain, then $R$ is not Teter type.
%\end{thm}

%\begin{lem}
%If $R$ is a non-Gorenstein domain,
%then $\omega_R \ncong \mm$.
%\begin{proof}
%We assume there exists
%an isomorphism
%$\phi : \omega_R \rightarrow \mm$.
%Then there exists a homomorphism
%$\iota : \omega_R \rightarrow R$
%such that $\iota(\omega_R)=\mm$.
%Thus $\mm \subseteq \tr(\omega_R)$.
%Since $R$ is not Gorenstein,
%we get $\tr(\omega_R)=\mm=\iota(\omega_R)$.
%Then $R$ is Teter type,
%and this contradicts to Theorem 2.4.
%\end{proof}
%\end{lem}

\begin{prop}\label{2.3}
Let $S$ be a Cohen-Macaulay homogeneous affine semigroup,
and let $G_S=\{ \bfa_1,\cdots,\bfa_s \} \subseteq \mathbb{N}^n$ be the minimal generators of $S$, where $s \geqq 2$.
If $R$ is non-Gorenstein and nearly Gorenstein,
then $2 \leqq |V_{\min}(S)| \leqq s-1$.
\begin{proof}
$|V_{\min}(S)| \geqq 2$ follows from by Theorem 2.2.
We show $|V_{\min}(S)| \leqq s-1$.
Assume $|V_{\min}(S)| \geqq s$.
Since $R$ is nearly Gorenstein,
there exist $\bfv_1 \in V_{\min}(S)$
and $\bfu_1 \in S-V(S)$
such that $\bfa_1=\bfu_1+\bfv_1$ by Proposition \ref{2.1}.
Since $\bfu_1 \in S-V(S)$,
we get $\bfu_1+\bfv \in S$ and
$\deg \bfx^{\bfu_1+\bfv}=1$ for all $\bfv \in V_{\min}(S)$.
Thus we obtain $|V_{\min}(S)|=s$ and
$\{\bfu_1+\bfv_i; 1 \leqq i \leqq s\}=G_S$.
Then we have $\omega_R \cong \mm$.
This contradicts to \cite[Chapter I, Theorem 12.9]{Sta}.
Therefore, we get $|V_{\min}(S)| \leqq s-1$.
  \end{proof}
  \end{prop}

\subsection{Numerical semigroup rings}
A semigroup $S =  \langle a_1 ,\cdots, a_n  \rangle$ with $a_i \in \NN$ is called \textit{numerical semigroup} if $\gcd(a_1,\cdots,a_n) = 1$.
We can assume $0<a_1 < a_2 < \cdots < a_n$.
We call $n$ \textit{embedding dimension} of $S$.
%The condition $\gcd(a_1,\cdots,a_n) = 1$
%is equivalent to say that
There exists $s \in S$ such that $s+S \subseteq S$.
The smallest such integer $c$ is called \textit{conductor} of $S$,
and we call $S$ \textit{symmetric} if the number of elements in $S$ that are less than $c$ is ${c}/2$.
The \textit{Frobenius number} of $S$ is denoted by $F(S)$
and it is the conductor of $S$ minus one.
Moreover, $\PF(S):=\{ x \notin S : x + s \in S \; \textit{for any}\; s \in S \}$
is called \textit{the set of pseudo-Frobenius numbers} of $S$ (see \cite{71}).

\begin{prop}[{see \cite[Corollary 4.11]{NM}}]\label{2.4}
Let $S$ be a numerical semigroup. The following are equivalent:
\begin{itemize}
\item[(1)] $S$ is symmetric;
\item[(2)] $\PF(S) = \{F(S)\}$;
\item[(3)] $\kk[S]$ is Gorenstein.
\end{itemize}
\end{prop}

\begin{prop}[{see \cite[Proposition 2.13 and Corollary 4.7]{NM}}]\label{2.5}
Every numerical semigroup of embedding dimension two $S=\langle a,b\rangle$ is symmetric.
In this case, $F(S)=ab-a-b$.
\end{prop}

%\begin{dfn}
%the Ap\'{e}ry set of $S$ with
%respect to $s \in S\setminus \{0\}$ is
%$\Ape(S,s) = \{x\in S :x-s \notin S\}$.
%\end{dfn}
\begin{dfn}[{see \cite[Definition 2.3]{pm1}}]
Let $S$ be a numerical semigroup and let $s \in S \setminus \{0\}$;
the \textit{Ap\'{e}ry} set with respect to $s$ of $S$ is the generating set
$B =
\{b_0 ,...,b_{s-1}\}$ such that $b_0=s$ and,
for $i>0$ , $b_i$ is the least
integer in $S$ having $s$-residue distinct from those of
$b_0 ,...,b_{i-1}$.
\end{dfn}

\begin{prop}[{see \cite[Lemma 3.1]{pm1}}]\label{2.7}
Let $S =  \langle a_1 ,...,a_n \rangle$ be a numerical semigroup.
Let $s \in S \setminus \{0\}$, and assume
$\sum_{i=1}^n c_i a_i \in \Ape(s,S)$.
Then $\sum_{i=1}^n d_i a_i \in \Ape(s,S)$ whenever
$0 \leqq {d_i} \leqq c_i$.
\end{prop}

\subsection{Projective monomial curves}
Let $\bfa=a_1,a_2,\cdots,a_n$ be a sequence
of positive integers with $\gcd(a_1,a_2,\cdots,a_n)=1$ and $a_1<a_2<\cdots<a_n$.
Then, we define a homogeneous affine semigroup
$$S_\bfa =  \langle(0,a_n),(a_1,a_n-a_1),(a_2,a_n-a_2),\cdots,
(a_{n-1},a_n-a_{n-1}),(a_n,0) \rangle.$$
We call $S_\bfa$
{the \textit{projective monomial curve}} defined by $\bfa$.
Moreover, for each $i=1,2$, we put
$S_i = \pi_i(S_{\bfa})$, where $\pi_i$ is the natural projection to the
$i$-th component.
Then $S_1$ and $S_2$ are numerical {semigroups}.

\begin{dfn}[{see \cite[Definitions 4.4]{pm1}}]
Let $B_i=\Ape(a_n,S_i)$ be the Ap\'{e}ry set with respect to $a_n$ of $S_i$ for each $i=1,2$.
We call \textit{Ap\'{e}ry set} with respect to $a_n$ of $S_{\bfa}$ the generating set
$\Ape(S_{\bfa}) = \{b_0,b_1,\cdots,b_{a_n}\}$
where $b_0 = (0,a_n),
b_{a_n} = (a_n,0)$
and $b_i = (\nu_i,\mu_i)$ satisfy the following for any $1 \leqq i \leqq a_n-1$:
\begin{itemize}
\item[(i)] $\{a_n,\nu_1,\cdots,\nu_{a_n-1}\} = B_1$;
\item[(ii)] $\mu_i$ is the least element of $S_2$ such that $(\nu_i,\mu_i) \in S$.
\end{itemize}
The Ap\'{e}ry set $\Ape(S_{\bfa})$ is called \textit{good} if
$\{a_n,\mu_1,\mu_2,\cdots,\mu_{a_{n-1}}\} = B_2$.
\end{dfn}

\begin{prop}[{see \cite[Lemma 4.6]{pm1}}]\label{2.9}
The following conditions are equivalent:
\begin{itemize}
\item[(1)] $S_{\bfa}$ is a Cohen-Macaulay;
\item[(2)] $\Ape(S_{\bfa})$ is good.
\end{itemize}
\end{prop}

We recall {that $r(S_\bfa)$ denotes the} Cohen-Macaulay type of $\kk[S_{\bfa}]$.
\begin{prop}[{see \cite[Theorem 4.9]{pm1}}]\label{2.10}
Assume that $\kk[S_\bfa]$ is Cohen-Macaulay. Let $B$ be the
Ap\'{e}ry set with respect to $a_n$ of $S_{\bfa}$. Then
$r(S_{\bfa}) = |\tilde{B}|,$
where
$$\tilde{B} = \{b \in B \setminus
\{(0,a_n),(a_n,0)\} ; b+x \notin B \;\textit{for all}\;\; {x} \in B\}.$$
\end{prop}

\begin{prop}[{see \cite[Proposition 4.11]{pm1}}]\label{2.11}
Assume that $S_{\bfa}$
is Cohen-Macaulay and let
$B=\{(0,a_n),(a_n,0)\}\cup\{b_i=(\nu_i,\mu_i);1\leqq i \leqq a_n-1\}$ be the Ap\'{e}ry set with respect to $a_n$ of $S_{\bfa}$
ordered so that
$\nu_1<\cdots<\nu_{a_n-1}$.
Then $S_{\bfa}$ is Gorenstein if and only if
$b_{a_{n-1}} = b_i + b_{a_n-1-i}$ for all $i=1,\cdots,n-2$.
\end{prop}

\begin{dfn}[{see \cite[Section 3]{Gre}}]
Let
$\bfa = a_1,a_2,\cdots,a_n$
be a sequence of positive integers with
$\gcd(a_1,a_2,\cdots,a_n) = 1$ and
$a_1<a_2<\cdots<a_n$.
We define the dual sequence
$\bfa' = a_n-a_{n-1}, a_n-a_{n-2},\cdots,a_n-a_1,a_n$.
It is known that
$S_\bfa$ and $S_\bfa'$ are isomorphic.
\end{dfn}

Put $F_1 = \NN(a_n,0)$,
$F_2 = \NN(0,a_n)$ and put
$C_i=\{w \in \ZZ S_\bfa \; ; \;w+g \notin S_\bfa \; \textit{for any}\; g \in F_i\}$
for $i=1,2$, respectively.
Denote by $\kk[\omega_{S_\bfa}]$ the $R$-submodule of $\kk[\ZZ S_\bfa]$
generated by $\{ \bfx^{v} ; v \in \omega_{S_\bfa} \}$,
where $\omega_S=-(C_1 \cap C_2)$.
By applying \cite{Aff1} to our case,
the following is true.

\begin{prop}[{\cite[Theorem 3.8]{Aff1}}]\label{2.13}
If $\kk[S_\bfa]$ is Cohen-Macaulay, then $\kk[\omega_{S_\bfa}]$ is the canonical module of $\kk[S_\bfa]$.
\end{prop}

\begin{prop}[{\cite[Theorem 2.6]{Aff1}}]\label{2.14}
The following conditions are equivalent:
\begin{itemize}
\item[(1)] $\kk[S_\bfa]$ is not Cohen-Macaulay;
\item[(2)] There exists $w \in \ZZ{S_{\bfa}}\setminus S_{\bfa}$
such that $w+(0,a_n) \in S_{\bfa}$ and $w+(a_n,0) \in S_{\bfa}$.
\end{itemize}
\end{prop}

\subsection{Nearly Gorenstein movement of projective monomial curves}
Let $S_\bfa$ be the projective monomial curve.
Then by using Proposition \ref{2.1},
for any $\bfa_i \in G(S_\bfa)$,
there exists $\bfu \in S_{\bfa} - V(S_{\bfa})$
such that $\bfa_i=V_{\min}(S_{\bfa})+\bfu$.
Here $V_{\min}(S_{\bfa})+\bfu=\{\bfv+\bfu:\bfv \in V_{\min}(S_{\bfa})\}$.
For all $\bfa_i {\in G(S_\bfa)}$,
there exists such a covering
$V_{\min}(S_{\bfa})+\bfu$
of minimal generators
of canonical module.
Based on this, we introduce Nearly Gorenstein movement to prove Theorem A.

\begin{dfn}
Let $U \subseteq S_{\bfa}-V(S_{\bfa})$,
and let $\pi_1$ be {the natural} projection $\pi_1: \NN^2 \rightarrow \NN$ defined by $\pi_1{(a,b)}=a$.
We call $\move=\pi_1(U)$ {a} \textit{{nearly} Gorenstein movement} of {$V(S_{\bfa})$} if
$G_{S_{\bfa}}=\bigcup_{\bfu \in U} (V_{\min}(S_{\bfa})+\bfu).$
If $\move=\pi_1(U)$ {is a} nearly Gorenstein movement of ${V(S_{\bfa})}$,
we define \textit{nearly Gorenstein covering} $C_\move$ as
{$\{\pi_1(V_{\min}(S_{\bfa})+\bfu) \;:\; \bfu \in U \}$}.
\end{dfn}

By Proposition \ref{2.1}, the following is true.
\begin{prop}\label{2.16}
Let $S_\bfa$ be the projective monomial curve.
$S_\bfa$ is nearly Gorenstein
if and only if
there exists nearly Gorenstein movement of $S_\bfa$.
\end{prop}

\begin{ex}\label{2.17}
Let $\bfa=1,2,3,4$.
We have $V(S_{\bfa})=V_{\min}(S_{\bfa})=\{(1,3),(2,2),(3,1)\}$
since $\omega_{S_\bfa}=\left <(1,3),(2,2),(3,1) \right>$.
%In this case,
{If we set $U=\{(-1,1),(0,0),(1,-1)\} \subseteq S_{\bfa}-V(S_{\bfa})$,
then
$\move=\pi(U)=\{-1,0,1\}$} is a nearly Gorenstein movement of {$V(S_\bfa)$}.
{Indeed,
since $V_{\min}(S_{\bfa})+(-1,1)=\{(0,4),(1,3),(2,2)\}$
and
$V_{\min}(S_{\bfa})+(1,-1)=\{(2,2),(3,1),(4,0)\}$,
$$\bigcup_{\bfu \in U} (V_{\min}(S_{\bfa})+\bfu)
=
V_{\min} \cup (V_{\min}(S_{\bfa})+(-1,1)) \cup (V_{\min}(S_{\bfa})+(1,-1))
=G_{S_{\bfa}}
.$$
}

In this case,
nearly Gorenstein covering is
$C_\move=\{\{0,1,2\},\{1,2,3\},\{2,3,4\}\}$
{and $S_\bfa$ is nearly Gorenstein
because there is a nearly Gorenstein movement $\move=\{-1,0,1\}$ of $V(S_\bfa)$.}
The following figure {represents} this covering.
\end{ex}

\begin{figure}[htbp]
  \begin{minipage}[b]{0.45\linewidth}
    \centering
\begin{tikzpicture}
\draw [help lines] (0,0) grid (4.4,4.4);%(0,0)から(3.5,2.5)までの"細線の方眼"
\coordinate (O) at (0,0);
%\node at (O) [left] {$O$}; %(0,0)に点Oを設定し，その左に$O$と記載する
\coordinate (A) at (0,4);
\node at (A) [left, text=red] {$\underline{0}$};
\coordinate (B) at (1,3);
\node at (B) [left, text=red] {$\underline{1}$};
\coordinate (C) at (2,2);
\node at (2,2) [left, text=red] {$\underline{2}$};
\coordinate (D) at (3,1);
\node at (3,1) [left] {$3$};
\coordinate (E) at (4,0);
\node at (E) [left] {$4$};
%\fill (O) circle (2pt) (A) circle (2pt) (B) circle (2pt); %点O,A,Bに黒丸をつける（1つの\fillで複数の点を対象にできる）

\draw [->] (O)--(A); %矢印OAを描く
\draw [->] (O)--(B); 
\draw [->] (O)--(C); 
\draw [->] (O)--(D);
\draw [->] (O)--(E); 
\end{tikzpicture}
    \subcaption{$\{0,1,2\} \in C_\move$}
    
  \end{minipage}
  \begin{minipage}[b]{0.45\linewidth}
    \centering
 \begin{tikzpicture}
\draw [help lines] (0,0) grid (4.4,4.4);%(0,0)から(3.5,2.5)までの"細線の方眼"
\coordinate (O) at (0,0);
%\node at (O) [left] {$O$}; %(0,0)に点Oを設定し，その左に$O$と記載する
\coordinate (A) at (0,4);
\node at (A) [left] {$0$};
\coordinate (B) at (1,3);
\node at (B) [left, text=red] {$\underline{1}$};
\coordinate (C) at (2,2);
\node at (2,2) [left, text=red] {$\underline{2}$};
\coordinate (D) at (3,1);
\node at (3,1) [left, text=red] {$\underline{3}$};
\coordinate (E) at (4,0);
\node at (E) [left] {$4$};
%\fill (O) circle (2pt) (A) circle (2pt) (B) circle (2pt); %点O,A,Bに黒丸をつける（1つの\fillで複数の点を対象にできる）

\draw [->] (O)--(A); %矢印OAを描く
\draw [->] (O)--(B); 
\draw [->] (O)--(C); 
\draw [->] (O)--(D);
\draw [->] (O)--(E); 
\end{tikzpicture}
     \subcaption{$\{1,2,3\} \in C_{\move}$}
  \end{minipage}
    \begin{minipage}[c]{0.45\linewidth}
    \centering
\begin{tikzpicture}
\draw [help lines] (0,0) grid (4.4,4.4);%(0,0)から(3.5,2.5)までの"細線の方眼"
\coordinate (O) at (0,0);
%\node at (O) [left] {$O$}; %(0,0)に点Oを設定し，その左に$O$と記載する
\coordinate (A) at (0,4);
\node at (A) [left] {$0$};
\coordinate (B) at (1,3);
\node at (B) [left] {$1$};
\coordinate (C) at (2,2);
\node at (2,2) [left, text=red] {$\underline{2}$};
\coordinate (D) at (3,1);
\node at (3,1) [left, text=red] {$\underline{3}$};
\coordinate (E) at (4,0);
\node at (E) [left, text=red] {$\underline{4}$};
%\fill (O) circle (2pt) (A) circle (2pt) (B) circle (2pt); %点O,A,Bに黒丸をつける（1つの\fillで複数の点を対象にできる）

\draw [->] (O)--(A); %矢印OAを描く
\draw [->] (O)--(B); 
\draw [->] (O)--(C); 
\draw [->] (O)--(D);
\draw [->] (O)--(E); 
\end{tikzpicture}
    \subcaption{$\{2,3,4\}$}
  \end{minipage}
  \caption{nearly Gorenstein {covering}}
\end{figure}
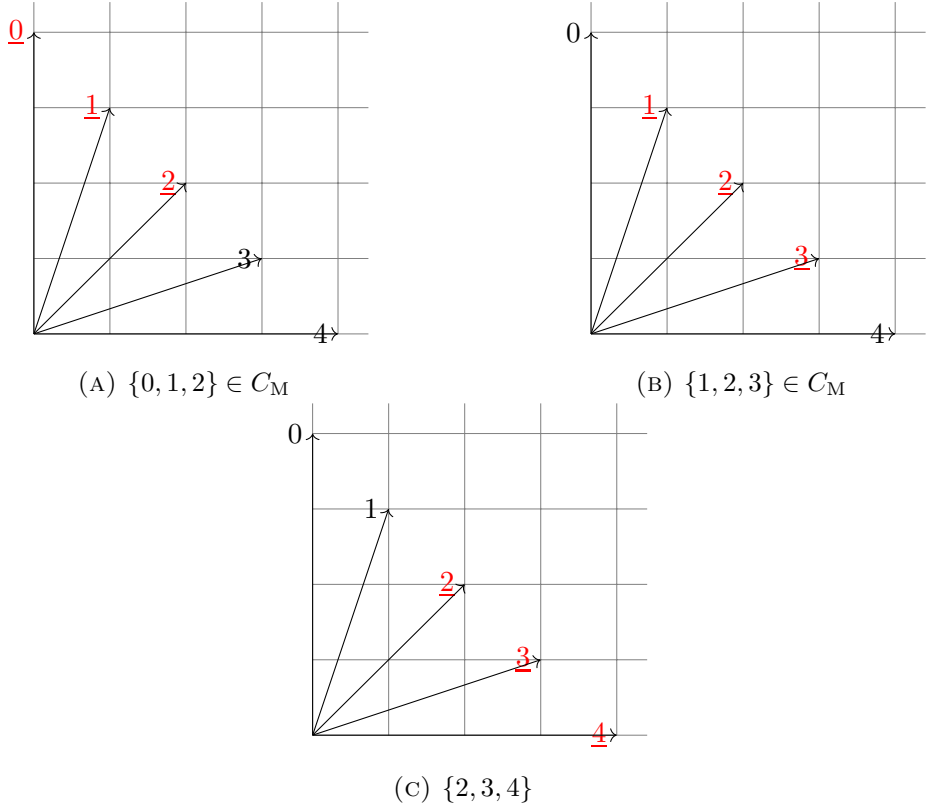

%{
%By Proposition 2.16, we know that for each $0 \leqq i \leqq 4$, there is a $C \in C_M$ such that $i \in C$.
%For example, $0,1,2 \in \{0,1,2\} \in C_M$ and
%$3,4 \in \{2,3,4\} \in C_M$.
%On the other hand,
%if we set $U'=\{(-1,1),(1,-1)\}$,
%then $\move'=\pi(U')=\{-1,1\}$ is also a nearly Gorenstein movement of $V(S_\bfa)$
%and
%its nearly Gorenstein covering is
%$C_{\move'}=\{\{0,1,2\},\{2,3,4\}\}$.}
%{If there exists 
%a nearly Gorenstein movement $\move$
%and its nearly Gorenstein covering $C_\move$,}
%we can simply write nearly Gorenstein covering as follows:

\begin{dfn}\label{2.18}
{Let $S_\bfa$ be the projective monomial curve defined by $\bfa=a_1,a_2,\cdots,a_n$ and let $C_\move$ be a nearly Gorenstein covering of $S_\bfa$.
For $X,X' \in C_\move$, define
$X \leqq X'$ as $\min X \leqq \min X'$.
Under this ordering, $C_\move$ is totally ordered set.
Thus we can write $C_\move=\{X_1,X_2,\cdots,X_m\}$ with $X_i<X_j$ for any $1 \leqq i < j \leqq m$.
If $X \leqq X'$, then $X'=\{x+\min X'-\min X \;:\; x \in X\}$.
Therefore, $C_\move$ can be written as follows:
$$X_1 \xrightarrow{\min X_2 - \min X_1} X_2 \xrightarrow {\min X_3 - \min X_2}
\cdots
\xrightarrow{\min X_{m-1} - \min X_{m-2}}
X_{m-1}
\xrightarrow{\min X_m - \min X_{m-1}} X_m \;$$
In this diagram, for each $1\leqq i \leqq m$, arrange natural numbers from $1$ to $n$ and represent $X_i$ by underlining only beneath the elements of $X_i$.
}
\end{dfn}

\begin{ex}
{Let $\bfa=1,2,3,4$. Then we observed that
$C_\move=\{\{0,1,2\},\{1,2,3\},\{2,3,4\}\}$
is a nearly Gorenstein covering in Example \ref{2.17}.
Following Definition \ref{2.18}, when represented graphically, $C_\move$ looks as follows:
$$\underline{0,1,2},3,4 \xrightarrow{1}
0,\underline{1,2,3},4 \xrightarrow{1}
0,1,\underline{2,3,4}$$}
%
%$\underline{0,1,2},3,4 \xrightarrow{2}
%0,1,\underline{2,3,4}$
\end{ex}

\section{Proof of Lemma A}
Now let us show the following lemma.

\begin{mainlem}\label{thm:lemmaA}
($\alpha$) Let $\bfa=k,k+1,2k+1$ for some $k\geqq 1$,
then $S_\bfa$ is non-Gorenstein and nearly Gorenstein.
Moreover, $S_{\bfa}$ is level with $r(S_{\bfa})=2$.

($\beta$) Let $\bfa=a_1,a_2,a_3,a_4$ be a sequence of
positive {integers} such that $\bfa \neq 1,2,3,4$
with $\gcd(a_1,a_2,a_3,a_4)=1$ and $a_1<a_2<a_3<a_4$.
Then the following is true.
\begin{itemize}
\item[(i)] When $\bfa=a,b,a+b,a+2b$
for some positive integers $a<b$.
\begin{itemize}
\item[(a)] If $(a,b)=(k,k+1)$ for some
$k \geqq 1$, then $S_\bfa$ is Gorenstein.
\item[(b)] If $(a,b)=(2k-1,2k+1)$ for some
$k \geqq 1$, then $S_\bfa$ is non-Gorenstein and nearly Gorenstein.
In this case, $S_{\bfa}$ is level with $r(S_{\bfa})=3$.
\item[(c)] In other cases, then $S_\bfa$ is not Cohen-Macaulay.
\end{itemize}
\item[(ii)] When $\bfa=a,b,a+b,2a+b$
for some positive integers $a<b$.
\begin{itemize}
\item[(d)] If $(a,b)=(k,k+1)$ ($k \geqq 2$),
then $S_\bfa$ is not nearly Gorenstein.
\item[(e)] In other cases, $S_\bfa$ is not Cohen-Macaulay.
\end{itemize}
\item[(iii)]
If $\bfa=a,b,a+b,2b$
for some positive integers $a<b$, then $S_\bfa$ is not Cohen-Macaulay.
\item[(iv)] When $\bfa=a,b,2a,a+b$
for some positive integers $a<b<2a$.
\begin{itemize}
\item[(a)$'$] If $(a,b)=(2k+1,4k+1), (2k,4k-1)$ for some $k \geqq 1$,
then $S_\bfa$ is Gorenstein.
\item[(b)$'$] If $(a,b)=(2k+1,4k)$ for some $k \geqq 1$,
then $S_\bfa$ is non-Gorenstein and nearly Gorenstein.
In this case, $S_{\bfa}$ is level with $r(S_{\bfa})=3$.
\item[(c)$'$] In other cases, $S_\bfa$ is not Cohen-Macaulay.
\end{itemize}
\item[(v)] When $\bfa=a,2a,b,a+b$
for some positive integers $a<2a<b$.
\begin{itemize}
\item[(d)$'$] If $(a,b)=(k,2k+1)$ for some $k \geqq 2$,
then $S_\bfa$ is Cohen-Macaulay but not nearly Gorenstein.
\item[(e)$'$] In other cases, $S_\bfa$ is not Cohen-Macaulay.
\end{itemize}
\end{itemize}
\begin{proof}
($\alpha$)
Put $\bfa=k,k+1,2k+1$ for some $k\geqq 1$.
We show
$$S_{\bfa}=\langle (0,2k+1),(k,k+1),(k+1,k),(2k+1,0)\rangle
$$
is non-Gorenstein and nearly Gorenstein.
If $k=1$, in the same way as Example \ref{2.17},
$S_{\bfa}$ is nearly Gorenstein and level with $r(S_{\bfa})=2$.
Then we can assume $k\geqq 2$.
In this case, $S_1=S_2=\langle k,k+1\rangle$.
Now we put $H=\langle k,k+1 \rangle$.

First, we show $B_H=\Ape(2k+1,H)$ is equal to the following set $X$, where
$$X=\{2k+1\}\cup\{nk : 1 \leqq n \leqq k \}\cup\{ n(k+1) : 1 \leqq n \leqq k\}.$$
Indeed, since $k+(k+1)=2k+1$,
for any $s \in H$, the remainder of $s$ divided by $2k+1$ can be written like
$nk$ or $n(k+1)$ for some $n \geqq 0$.
Note that $(k+2)k=(k-1)(k+1)+(2k+1) \notin B_H$
and $(k+1)(k+1)=k^2+(2k+1) \notin B_H$.
So, by Proposition \ref{2.7}, we obtain the following:
\begin{itemize}
\item $nk \in B_H$ implies $0 \leqq n \leqq k$ or $n=k(k+1)$;
\item $n(k+1) \in B_H$ implies $n \leqq k$.
\end{itemize}
Then $B_H \subseteq X$ and since $|B_H|=|X|=2k+1$, we get $B_H=X$.
Therefore, we get $\Ape(2k+1,S_1)=\Ape(2k+1,S_2)=X$.

Next, we show $B_{S_\bfa}=\Ape(2k+1,S_\bfa)$ is equal to the following set $Y$, where
$$Y=\{(2k+1,0),(0,2k+1)\}\cup\{n(k,k+1);1 \leqq n \leqq k\}\cup\{ n(k+1,k) ; 1 \leqq n \leqq k \}.$$
It is enough to check that $n(k+1)-m(2k+1) \notin S_2$ and
$nk - m(2k+1) \notin S_2$ for any $1 \leqq n \leqq k-1$ and $m \geqq 1$.
Indeed, if $n(k+1)-m(2k+1) \in S_2$,
we can write $n(k+1)-m(2k+1)=ak+b(k+1)$ for some $a,b \in \NN$.
Then $(n-b-m)(k+1)=(a+m)k \neq 0$.
Since $\gcd(k,k+1)=1$, we can write $n-b-m=kl$ for some $0<l\in \NN$.
Thus $n=kl+b+m>k$, as a contradiction.
Therefore, $n(k+1)-m(2k+1) \notin S_2$ for any $1 \leqq n \leqq k-1$ and $m \geqq 1$.
By the same discussion, we get
$nk - m(2k+1) \notin S_2$ for any $1 \leqq n \leqq k-1$ and $m \geqq 1$.
Then we obtain $B_{S_\bfa}=Y$.
Since $\Ape(2k+1,S_2)=\pi_2(B_{S_{\bfa}})$, thus $S_\bfa$ is good,
so $S_\bfa$ is Cohen-Macaulay by using Proposition \ref{2.9}.
Now we show $S_\bfa$ is non-Gorenstein and nearly Gorenstein.
It is easy to check
$${\tilde{B}}_{S_\bfa}=\{k(k,k+1),k(k+1,k)\}.$$
Indeed,
it is true if $k=1$. We consider the case of $k \geqq 2$.
Since every degree of element of $B$ is less than or equal to $k$, we obtain $k(k,k+1),k(k+1,k) \in {\tilde{B}}_{S_\bfa}$.
On the other hand, the following is true any $1 \leqq n \leqq k-1$.
$$n(k,k+1)+(k,k+1)=(n+1)(k,k+1) \in {\tilde{B}}_{S_\bfa},$$
$$n(k+1,k)+(k+1,k)=(n+1)(k+1,k) \in {\tilde{B}}_{S_\bfa}.$$
Thus ${\tilde{B}}_{S_\bfa}=\{k(k,k+1),k(k+1,k)\}.$
Therefore, by Proposition \ref{2.10},
we get $r({S_\bfa})=2$.
%Since $\deg \bfx^{k(k,k+1)}=\deg \bfx^{k(k+1,k)}=k < \deg \bfx^{\bfb}$ for an
%$n(k,k+1)+(k,k+1)=(n+1)(k,k+1) \in B$ and
%$n(k+1,k)+(k+1,k)=(n+1)(k+1,k)$ for any $1 \leqq n \leqq k-1$,
Recall $\omega_S=-(C_1\cap C_2)$, 
where $F_1=\NN(2k+1,0)$, $F_2=\NN(0,2k+1)$
and $C_i=\{w \in \ZZ S_\bfa \;;\; w+g \notin S_\bfa \; \textit{for any}\; g \in F_i \}$ for $i=1,2$
(see Proposition \ref{2.13}).
Next, we show
$$\omega_{S_\bfa}=\langle((-(k^2-k-1),-(k^2-2k-1)),
(-(k^2-2k-1),-(k^2-k-1))\rangle.$$
Indeed, by Proposition \ref{2.5},
$F(H)=k^2-k-1$.
Thus
$$C_1\cap C_2 \subseteq \{ (x,y) \in \ZZ{S_{\bfa}} : x \leqq k^2-k-1 \; \textit{and}\; y \leqq k^2-k-1 \}.$$
Therefore,
\begin{equation}
\omega_{S_\bfa} \subseteq
Z:=\{ (x,y) \in \ZZ{S_{\bfa}} : x \geqq -(k^2-k-1) \; \textit{and}\; y \geqq -(k^2-k-1) \}.
\end{equation}
On the other hand, since
$k^2-k-1 = F(H) \notin S_1 \cup S_2$ and
$k^2-2k-1= F(H)-k \notin S_1 \cup S_2$,
we obtain
\begin{equation}
(-(k^2-k-1),-(k^2-2k-1)),(-(k^2-2k-1),-(k^2-k-1)) \in \omega_{S_\bfa} \cap Z.
\end{equation}
From (1) and (2),
$(-(k^2-k-1),-(k^2-2k-1))$ and $(-(k^2-2k-1),-(k^2-k-1))$
{belong} to the system of minimal generators of $\omega_{S_{\bfa}}$.
Moreover, since $r(S_{\bfa})=2$, we get
$$\omega_{S_\bfa}=\langle(-(k^2-k-1),-(k^2-2k-1)),(-(k^2-2k-1),-(k^2-k-1))\rangle.$$
Then by Proposition \ref{2.1}, $S_{\bfa}$ is non-Gorenstein and nearly Gorenstein.
In particular, since every element of $\omega_S$ has degree $-k+1$,
$S_\bfa$ is level with $r(S_\bfa)=2$.
%==============================================================
\bigskip

($\beta$)
(i) (a)
Put $\bfa=k,k+1,2k+1,3k+2$ for some $k\geqq 1$.
We show $$S_{\bfa}=\langle(0,3k+2),(k,2k+2),(k+1,2k+1),(2k+1,k+1),(3k+2,0)\rangle$$
is Gorenstein.
In this case, $S_1= \langle k,k+1 \rangle$ and
$S_2= \langle k+1,2k+1 \rangle$.
\begin{itemize}
\item First, we show $B_{S_1}=\Ape(3k+2,S_1)$ is equal to the following set $X$, where
$$X=\{3k+2\}\cup\{nk : 1 \leqq n \leqq k \}\cup\{ n(k+1) : 1 \leqq n \leqq k+1\} \cup \{ nk+k+1 : 1 \leqq n \leqq k \}.$$
%=================================================================
It is easy to check if $k=1$, then now we assume $k \geqq 2$.
Since $k+2(k+1)=3k+2$,
for any $s \in S_1$, the remainder of $s$ divided by $3k+2$ can be written like
$nk$ or $n(k+1)$ or $nk+k+1$ for some $n \geqq 0$.
Note that $(k+2)k=(k-2)(k+1)+(3k+2) \notin B_{S_1}$,
$(k+2)(k+1)=k^2+(3k+2) \notin B_{S_1}$
and
$(k+2)k+(k+1)=(k-1)(k+1)+3k+2 \notin B_{S_1}$.
So, in the same way as ($\alpha$),
%by Proposition 2.10, we obtain the following:
%\begin{itemize}
%\item[(1)] $nk \in B_H$ implies $0 \leqq n \leqq k$ or $n=k(k+1)$;
%\item[(2)] $n(k+1) \in B_H$ implies $0 \leqq n \leqq k+1$;
%\item[(3)] $nk+(k+1) \in B_H$ implies $0 \leqq n \leqq k$ or $n=(k+1)^2$.
%\end{itemize}
we get $B_{S_1} \subseteq X$.
Moreover, it is easy to check that $|B_{S_1}|=|X|=3k+2$, thus we get $B_{S_1}=X$.
\end{itemize}
%============================================================
In the same way as above, we obtain
$B_{S_2}=\Ape(3k+2,S_2)$ is equal to the following set $Y$, where
$$Y=\{3k+2\}\cup\{n(k+1) : 1 \leqq n \leqq 2k \}\cup\{ n(2k+1) : 1 \leqq n \leqq k+1\}.$$
%===============================================================
Next, we show the Ap\'{e}ry set $B_{S_\bfa}$ with respect to $3k+2$ of $S_\bfa$ is
\begin{equation*}
\begin{split}
Z= \{(3k+2,0),(0,3k+2)\} &\cup\{n(k,2k+2) : 1 \leqq n \leqq k \} \cup 
\{ n(k+1,2k+1) : 1 \leqq n \leqq k+1\} \\
& \cup \{ (n-1)(k,2k+2)+(2k+1,k+1) : 1 \leqq n \leqq k+1\}.
\end{split}
\end{equation*}
It is easy to check if $k=1$.
If $k \geqq 2$,
in the same way as ($\alpha$),
we can check that $n(2k+2)-m(3k+2) \notin S_2$,
$n'(2k+1) - m(3k+2) \notin S_2$
and $(n'-1)(2k+2)+k+1-m(3k+2) \notin S_2$
for any $1 \leqq n \leqq k-1$,
$1 \leqq n' \leqq k$ and $m \geqq 1$.
%Indeed, if $n(k+1)-m(2k+1) \in S_2$,
%we can write $n(k+1)-m(2k+1)=ak+b(k+1)$ for some $a,b \in \NN$.
%Then $(n-b-m)(k+1)=(a+m)k \neq 0$.
%Since $\gcd(k,k+1)=1$, we can write $n-b-m=kl$ for some $0<l\in %\NN$.
%Thus $n=kl+b+m>k$, as a contradiction.
%Therefore, $n(k+1)-m(2k+1) \notin S_2$ for any $1 \leqq n \leqq %k-1$ and $m \geqq 1$.
%By the same discussion, we get
%$nk - m(2k+1) \notin S_2$ for any $1 \leqq n \leqq k-1$ and $m %\geqq 1$.
Then we obtain $B_{S_\bfa}=Z$.
Thus $S_\bfa$ is good, so $S_\bfa$ is Cohen-Macaulay by Proposition \ref{2.9}.
Moreover,
\begin{equation}
\begin{split}\nonumber
(k+1)(k+1,2k+1)&=n(k,2k+2)+(k-n)(k,2k+2)+(2k+1,k+1) \\
&=n(k+1,2k+1)+(k+1-n)(k+1,2k+1)
\end{split}
\end{equation}
for any $1 \leqq n \leqq k$.
Therefore, $S_\bfa$ is Gorenstein by Proposition \ref{2.11}.

(b)
Put $\bfa=2k-1,2k+1,4k,6k+1$ for some $k\geqq 1$.
We show $$S_{\bfa}=\langle(0,6k+1),(2k-1,4k+2),(2k+1,4k),(4k,2k+1),(6k+1,0)\rangle$$
is non-Gorenstein and nearly Gorenstein.
By the same discussion as above,
we get the following.
\begin{equation*}
\begin{split}
\Ape(6k+1,S_1)=\{6k+1\}&\cup\{n(2k-1) : 1 \leqq n \leqq 2k \}\cup\{ n(2k+1) : 1 \leqq n \leqq 2k\} \\
&\cup \{ n(2k-1)+2k+1 : 1 \leqq n \leqq 2k \}, 
\end{split}
\end{equation*}
$$\Ape(6k+1,S_2)=\{6k+1\}\cup\{n(2k+1) : 1 \leqq n \leqq 4k \}\cup\{ n(4k) : 1 \leqq n \leqq 2k\},$$
\begin{equation*}
\begin{split}
B=\{(6k+1,0),(0,6k+1)\}&\cup\{n(2k-1,4k+2) : 1 \leqq n \leqq 2k \}\cup\{ n(2k+1,4k) : 1 \leqq n \leqq 2k\} \\
&\cup\{ (n-1)(2k-1,4k+2)+(4k,2k+1) : 1 \leqq n \leqq 2k\}.
\end{split}
\end{equation*}
Here $S_1= \langle2k-1,2k+1 \rangle$,
$S_2= \langle2k+1,4k \rangle$
and
$B=\Ape(6k+1,S_\bfa)$.
Thus $S_{\bfa}$ is Cohen-Macaulay by Proposition \ref{2.9}.
Now we show $S_{\bfa}$ is non-Gorenstein and nearly Gorenstein.
In the same way as ($\alpha$),
it is easy to check
$$\tilde{B}=\{2k(2k-1,4k+2),2k(2k+1,4k),(2k-1)(2k-1,4k+2)+(4k,2k+1)\}.$$
Therefore,
we get $r({S_\bfa})=3$ by Proposition \ref{2.10}.
Next we show
$\omega_{S_\bfa}= \langle
v_1,
v_2,
v_3
\rangle,$
where
$v_1=(-(4k^2-4k-1),-(8k^2-6k-1))$,
$v_2=(-(4k^2-6k),-(8k^2-4k-2))$
and
$v_3=(-(4k^2-8k-1),-(8k^2-2k-1))$.
Indeed, by Proposition \ref{2.5},
$F(S_1)=4k^2-4k-1$ and $F(S_2)=8k^2-2k-1$.
Thus we obtain the following.
\begin{equation}
\omega_{S_\bfa} \subseteq
X:=\{ (x,y) \in \ZZ{S_{\bfa}} : x \geqq -(4k^2-4k-1) \; \textit{and}\; y \geqq -(8k^2-2k-1) \}.
\end{equation}
On the other hand, the following is true.
$$-\pi_1(v_1) = F(S_1) \notin S_1,
-\pi_1(v_2) = F(S_1)-(2k-1) \notin S_1,
-\pi_1(v_3) = F(S_1) - 4k \notin S_1,$$
$$-\pi_2(v_1)= F(S_2) - 4k \notin S_2,
-\pi_2(v_2) = F(S_2)-(2k+1) \notin S_2,
-\pi_2(v_3) = F(S_2) \notin S_2.$$
Therefore, we get
$
v_1,v_2,v_3
\in \omega_{S_\bfa} \cap Z$.
From this and (3),
$v_1,v_2$ and $v_3$
are belong to the system of minimal generators of $\omega_{S_{\bfa}}$.
Moreover, since $r(S_{\bfa})=3$, we get
$\omega_{S_\bfa}=\langle
v_1,v_2,v_3
\rangle.$
Then by Proposition \ref{2.1}, $S_{\bfa}$ is non-Gorenstein and nearly Gorenstein.
In particular, since every element of $\omega_S$ has degree $-2k+2$,
$S_\bfa$ is level with $r(S_\bfa)=3$.

(c) Put $\bfa=a,b,a+b,a+2b$ where $b \geqq a+3$ and $\gcd(a,b)=1$.
We show
$$S_{\bfa}=\langle(0,a+2b),(a,2b),(b,a+b),(a+b,b),(a+2b,0) \rangle$$
is not Cohen-Macaulay.
Put $v=(a(b-1),2b^2-a-4b)$.
Since $v+(0,a+2b)=(b-1)(a,2b)$ and
$v+(a+2b,0)=(b-a-3)(0,a+2b)+(a+2)(b,a+b)$,
it is enough to show $v \notin S_\bfa$ by Proposition \ref{2.14}.
Assume that $v \in S_\bfa$.
Then {there exist} $c_1,\cdots,c_5 \in \NN$
such that $\sum_{i=1}^5c_if_i=v$.
Since $\pi_1(\sum_{i=1}^5c_if_i)=\pi_1(v)$, we get
\begin{equation}
(c_2+c_3+2c_4+3c_5)a+(c_3+c_4+2c_5)b=(b-1)a.
\end{equation}
Since $\gcd(a,b)=1$,
$c_3+c_4+2c_5 \equiv 0 \pmod a$.
Here, $c_3+c_4+2c_5 \neq 0$.
Indeed, if $c_3+c_4+2c_5 = 0$,
we get $c_3=c_4=c_5=0$.
Then $c_1f_1+c_2v_2=v$.
However, since $\deg \bfx^v= b-2$,
$\pi_1(c_1f_1+c_2v_2)=c_2a \leqq (b-2)a < \pi_1(v)$.
This yields a contradiction.
Thus we can write
$c_3+c_4+2c_5=al$, where $0<l \in \NN$, and
substitute it into (4),
we get
$c_2+c_4+c_5+al + bl =b-1$.
Then
$0 \leqq c_2+c_4+c_5 + l(a-1) + b(l-1) = -1 <0$, as a contradiction.

(ii)
(d)
Put $\bfa=k,k+1,2k+1,3k+1$ for some $k\geqq 1$.
We show
$$S_{\bfa}=\langle(0,3k+1),(k,2k+1),(k+1,2k),(2k+1,k),(3k+1,0)\rangle$$ is not nearly Gorenstein.
In the same way as above,
we get the following.
$$\Ape(3k+1,S_1)=\{3k+1\}\cup\{nk : 1 \leqq n \leqq k \}\cup\{ n(k+1) : 1 \leqq n \leqq k\} \cup \{ k+n(k+1) : 1 \leqq n \leqq k \},$$
$$\Ape(3k+1,S_2)=\{3k+1\}\cup\{nk : 1 \leqq n \leqq 2k \}\cup\{ n(2k+1) : 1 \leqq n \leqq k\},$$
\begin{equation*}
\begin{split}
B=\{(3k+1,0),(0,3k+1)\}&\cup\{n(k,2k+1) : 1 \leqq n \leqq k \}
\cup\{ n(k+1,2k) : 1 \leqq n \leqq k\} \\
&\cup\{ (n-1)(k+1,2k)+(2k+1,k) : 1 \leqq n \leqq k\}.
\end{split}
\end{equation*}
Here $S_1= \langle k,k+1 \rangle,
S_2= \langle k,2k \rangle$
and $B=\Ape(3k+1,S_{\bfa})$.
Thus $S_{\bfa}$ is Cohen-Macaulay by Proposition \ref{2.9}.
In the same way as (b),
we can check
$$\tilde{B}=\{k(k,2k+1),k(k+1,2k),(k-1)(k+1,2k)+(2k+1,k)\},\;
\omega_{S_\bfa}= \langle v_1, v_2, v_3\rangle$$
where
$v_1=(-(4k^2-4k-1),-(8k^2-6k-1)),
v_2=(-(4k^2-6k),-(8k^2-4k-2))$
and
$v_3=((-(4k^2-8k-1),-(8k^2-2k-1))$.
Now we assume $S_{\bfa}$ is nearly Gorenstein.
Then by Proposition \ref{2.1},
there exist $u \in S_{\bfa}-V(S_{\bfa})$
and $v \in \{v_1,v_2,v_3\}$
such that
$v+u=(0,3k+1)$ and
$v+u \in S_{\bfa}$.
It is easy to check this yields a contradiction,
so $S_{\bfa}$ is not nearly Gorenstein.

(e)
Put $\bfa=a,b,a+b,2a+b$ where $b \geqq a+2$ and $\gcd(a,b)=1$.
Then
$$S_{\bfa}=\langle(0,2a+b),(a,a+b),(b,2a),(a+b,a),(2a+b,0)\rangle$$
is not Cohen-Macaulay.
Put $v=(a(b-1),b^2-2b+ab-3a)$.
Thus $v+(0,2a+b)=(b-1)(a,a+b)$ and
$v+(2a+b,0)=(b-a-2)(0,2a+b)+a(b,2a)+(a+b,a)$.
In the same way as (c),
we can check $v \notin S_\bfa$.
Then $S_{\bfa}$ is not Cohen-Macaulay by Proposition \ref{2.14}.
%Assume that $v \in S$.
%Then we there exists $c_1,\cdots,c_5 \in \NN$
%such that $\sum_{i=1}^5c_if_i=v$.
% Since $\pi_1(\sum_{i=1}^5c_if_i)=\pi_1(v)$, we get
% \begin{equation}
% (c_2+c_4+2c_5)a+(c_3+c_4+c_5)b=(b-1)a.
% \end{equation}
% Since $\gcd(a,b)=1$,
% $c_3+c_4+c_5 \equiv 0 \pmod a$.
% Here, $c_3+c_4+c_5 \neq 0$.
% Indeed, if $c_3+c_4+c_5 = 0$,
% we get $c_3=c_4=c_5=0$.
% Then $c_1f_1+c_2v_2=v$.
% However, since $\deg \bfx^v= b-2$,
% $\pi_1(c_1f_1+c_2v_2)=c_2a \leqq (b-2)a < \pi_1(v)$.
% This yields a contradiction.
% Thus we can write
% $c_3+c_4+c_5=al$, where $0<l \in \NN$, and
% substitute it for (2),
% we get
% $c_2+c_4+2c_5+bl = b-1$.
% Then
% $0 \leqq c_2+c_4+2c_5 + b(l-1) = -1 <0$, as a contradiction.

(iii)
 Put $\bfa=a,b,a+b,2b$ where $0<a<b$, $\gcd(a,b)=1$ and $b \neq 2$.
We show $$S_{\bfa}= \langle(0,2b),(a,2b-a),(b,b),(a+b,b-a),(2b,0) \rangle$$ is not Cohen-Macaulay.
Put $v=(2a,2b-2a)$.
Then $v+(0,2b)=2(a,2b-a)$ and
$v+(2b,0)=2(a+b,b-a)$.
In the same way as (c),
we can check $v \notin S_\bfa$.
Thus $S_{\bfa}$ is not Cohen-Macaulay by Proposition \ref{2.14}.
% it is enough to show $v \notin S$ by Proposition 2.17.
% Assume that $v \in S$.
% Since $\deg \bfx^v = 2b$,
% we get $v=f_i$ for some $1 \leqq i \leqq 5$.
% It is easy to check all case of this contradict to assumpition of $a,b$.

(iv)
(a)$'$
Put $\bfa=a,b,2a,a+b$
for some positive integers $a<b<2a$.
Then we get
$\bfa'=b_1,b_2,b_1+b_2,b_1+2b_2$
with $b_1<b_2$ and $\gcd(b_1,b_2)=1$,
here $b_1=b-a,b_2=a$.
Since $S_\bfa \cong S_{{\bfa}'}$,
we get (a)$'$ fron (a).
By using such a duality,
(b)$'$,(c)$'$,(d)$'$ and (e)$'$ are also followed
from (b),(c),(d) and (e), respectively.
  \end{proof}
  \end{mainlem}

\begin{ex}
(1) If $\bfa=6,7,13$,
then $S_{\bfa}$ is non-Gorenstein and nearly Gorenstein
by ($\alpha$) of Lemma A.
Since $\omega_{S_{\bfa}}=\langle (-29,-23),(-23,-29)\rangle$,
the nearly Gorenstein {covering} of $S_{\bfa}$ is as follows:
$\underline{0,6},7,13 \xrightarrow{7}
0,6,\underline{7,13}$

(2) If $\bfa=6,7,13,20$,
then $S_{\bfa}$ is Gorenstein
by (a) of Lemma A.
The nearly Gorenstein {covering} of $S_{\bfa}$ is as follows:
$\underline{0},7,13,20 \xrightarrow{7}
0,\underline{7},13,20 \xrightarrow{6}
0,7,\underline{13},20 \xrightarrow{7}
0,7,13,\underline{20}
$

(3) If $\bfa=5,7,12,19$,
then $S_{\bfa}$ is non-Gorenstein
and nearly Gorenstein
by (b) of Lemma A.
Since $\omega_{S_{\bfa}}=\langle (-23,-53),(-18,-58),(-11,-65) \rangle$,
the nearly Gorenstein {covering} of $S_{\bfa}$ is as follows:
$\underline{0,5},7,\underline{12},19 \xrightarrow{7}
0,5,\underline{7,12,19}
$

(4) If $\bfa=7,12,14,19$,
then $S_{\bfa}$ is non-Gorenstein
and nearly Gorenstein
by (b)$'$ of Lemma A.
Indeed, since $\bfa'=5,7,12,19$,
we get $S_{\bfa} \cong S_{\bfa'}$.
Thus $S_{\bfa}$ is non-Gorenstein
and nearly Gorenstein from (3).
\end{ex}
  
\section{Proof of Theorem A}
 Let $\bfa=a_1,\cdots,a_n$ be a sequence
 of positive integers with
 $\gcd(a_1,\cdots,a_n)=1$
 and $a_1<a_2<\cdots<a_n$.
First we consider the case of {$n=3$}.
{The following lemma appears in several papers. For example, refer to \cite[Lemma 4 and Section 4, Example 7]{BSV}, \cite[Theorem 3.3]{HQS} or \cite[Corollary 4.8]{dim2}.}

\begin{lem}\label{4.1}
  Let $\bfa=a,b,a+b$,
  then $S_{\bfa}$ is Cohen-Macaulay if and only if
  $b=a+1$.
\end{lem}

  \begin{prop}\label{4.2}
  Let $\bfa=a_1,a_2,a_3$.
  Then the following conditions are equivalent:
\item[(1)] $S_\bfa$ is non-Gorenstein and nearly Gorenstein;
\item[(2)] $\bfa=k,k+1,2k+1$ for some $k \geqq 1$.

If this is the case, $S_{\bfa}$ is level with $r(S_{\bfa})=2$.
  \begin{proof}
We know (2) implies (1) and $S_{\bfa}$ is level with $r(S_{\bfa})=2$ from ($\alpha$) of Lemma A.
Now we show (1) implies (2).
Assume $S_\bfa$ is non-Gorenstein and nearly Gorenstein.
  Then $|V_{\min}(S_{\bfa})|=2,3$ by Proposition \ref{2.3}.
When $|V_{\min}(S_{\bfa})|=3$.
{Since $S_\bfa$ is nearly Gorenstein, there exists nearly Gorenstein covering of $S_\bfa$ as follows.}
$${\underline{0,a_1,a_2},a_3 \xrightarrow{a_1}
0,\underline{a_1,a_2=2a_1,a_3=3a_1}}$$

{Since $\gcd(a_1,a_2,a_3)=1$, we get $\bfa=(a_1,a_2,a_3)=(a_1,2a_1,3a_1)=(1,2,3)$.}

  \item When $|V_{\min}(S_{\bfa})|=2$.
  Since $S_\bfa$ is nearly Gorenstein,
  there exists
  nearly Gorenstein {covering} of $S_\bfa$ as one of the following.
  \begin{itemize}
\item[(a)] $\underline{0,a_1},a_2,a_3 \xrightarrow{a_1}
0,\underline{a_1,a_2=2a_1},a_3 \xrightarrow{a_1}
0,a_1,\underline{2a_1,3a_1}
$
\item[(b)] $\underline{0,a_1},a_2,a_3 \xrightarrow{a_2}
0,a_1\underline{a_2,a_3=a_1+a_2}
$
\item[(c)] $\underline{0},a_1,\underline{a_2},a_3 \xrightarrow{a_1}
0,\underline{a_1},a_2,\underline{a_3=a_1+a_2}
$
\end{itemize}
In the case of (a), we get $\bfa=1,2,3$.
In the case of (b) and (c), we get $\bfa=a_1,a_2,a_1+a_2$.
Since $S_\bfa$ is Cohen-Macaulay,
we get $\bfa=(k,k+1,2k+1)$
by Lemma \ref{4.1}.
Therefore, we conclude that (1) implies (2).
\end{proof}
\end{prop}
Lastly, we consider the case of $n=4$.
\begin{prop}\label{4.3}
  Let $\bfa=a_1,a_2,a_3,a_4$.
  Then the following conditions are equivalent: 
\item[(1)] $S_\bfa$ is non-Gorenstein and nearly Gorenstein;
\item[(2)] $\bfa=1,2,3,4$ or $S_\bfa \cong S_{2k-1,2k+1,4k,6k+1}$ for some $k \geqq 1$.

If this is the case, $S_{\bfa}$ is level with $r(S_{\bfa})=3$.
  \begin{proof}
We know (2) implies (1) and $S_{\bfa}$ is level with $r(S_{\bfa})=3$ from (b) of Lemma A.
Now we show (1) implies (2).  
  Assume $S_\bfa$ is non-Gorenstein and nearly Gorenstein.
  Then $|V_{\min}(S_{\bfa})|=2,3,4$ by Proposition \ref{2.3}.
When $|V_{\min}(S_{\bfa})|=4$,
  in the same way as Proposition \ref{4.2},
we get
  $(a_1,a_2,a_3,a_4)=(1,2,3,4)$.
When $|V_{\min}(S_{\bfa})|=3$,
  since $S_\bfa$ is nearly Gorenstein, there exists nearly Gorenstein {covering} of $S_\bfa$ as one of the following.
\begin{itemize}
\item[(a)] $\underline{0,a_1,a_2},a_3,a_4 \xrightarrow{a_1}
0,\underline{a_1,a_2=2a_1,a_3=3a_1},a_4 \xrightarrow{a_2}
0,a_1,\underline{2a_1,3a_1,4a_1}
$
\item[(b)] $\underline{0,a_1,a_2},a_3,a_4 \xrightarrow{a_1}
0,\underline{a_1},a_2,\underline{a_3=2a_1,a_4=a_1+a_2}
$
\item[(c)] $\underline{0,a_1,a_2},a_3,a_4 \xrightarrow{a_2}
0,a_1,\underline{a_2,a_3=a_1+a_2,a_4=2a_2}
$
\item[(d)] $\underline{0,a_1},a_2,\underline{a_3},a_4 \xrightarrow{a_1}
0,\underline{a_1,a_2=2a_1},a_3,\underline{a_4=a_1+a_3}
$
\item[(e)] $\underline{0},a_1,\underline{a_2,a_3},a_4 \xrightarrow{a_1}
0,\underline{a_1},a_2,\underline{a_3=a_1+a_2,a_4=2a_1+a_2}
$
\end{itemize}
In the case of (a), we get $\bfa=1,2,3,4$.
In the case of (b), we get $\bfa=a,b,2a,a+b$ where $(a,b)=(a_1,a_2)$.
Thus we obtain $(a,b)=(2k+1,4k)$ for some
$k>0$ from (v) of Lemma A.
Then $S_{\bfa} \cong S_{\bfa'}=S_{2k-1,2k+1,4k,6k+1}$.
In the case of (c), we get
$\bfa=(a,b,a+b,2b)$ where $(a,b)=(a_1,a_2)$.
Thus we obtain $\bfa=1,2,3,4$ from (iii) of Lemma A.
In the case of (d), we get
$\bfa=(a,2a,b,a+b)$ where $(a,b)=(a_1,a_3)$.
Then we obtain $\bfa=1,2,3,4$ from (v) of Lemma A.
In the case of (e), we get
$\bfa=(a,b,a+b,a+2b)$ where $(a,b)=(a_1,a_2)$.
Thus we obtain $(a,b)=(2k-1,2k+1)$ for some
$k>0$ from (i) of Lemma A.
Then $S_{\bfa}=S_{2k-1,2k+1,4k,6k+1}$.
\item When $|V_{\min}(S_{\bfa})|=2$,
  since $S_\bfa$ is nearly Gorenstein, there exists nearly Gorenstein {covering} of $S_\bfa$ as one of the following.
\begin{itemize}
\item[(a)] $\underline{0,a_1},a_2,a_3,a_4 \xrightarrow{a_1}
0,\underline{a_1,a_2=2a_1},a_3,a_4 \xrightarrow{a_1}
0,a_1,\underline{2a_1,a_3=3a_1},a_4 \xrightarrow{a_1}
0,a_1,2a_1,\underline{3a_1,4a_1}
$
\item[(b)] $\underline{0,a_1},a_2,a_3,a_4 \xrightarrow{a_1}
0,\underline{a_1,a_2=2a_1},a_3,a_4 \xrightarrow{a_3-a_1}
0,a_1,2a_1,\underline{a_3,a_4=a_1+a_3}
$
\item[(c)] $\underline{0,a_1},a_2,a_3,a_4 \xrightarrow{a_2}
0,a_1,\underline{a_2,a_3=a_1+a_2},a_4 \xrightarrow{a_1}
0,a_1,a_2,\underline{a_1+a_2,a_4=2a_1+a_2}
$
\item[(d)] $\underline{0},a_1,\underline{a_2},a_3,a_4 \xrightarrow{a_1}
0,\underline{a_1},a_2,\underline{a_3=a_1+a_2},a_4 \xrightarrow{a_2-a_1}
0,a_1,\underline{a_2},a_1+a_2,\underline{a_4=2a_2}
$
\item[(e)] $\underline{0},a_1,\underline{a_2},a_3,a_4 \xrightarrow{a_1}
0,\underline{a_1},a_2,\underline{a_3=a_1+a_2},a_4 \xrightarrow{a_2}
0,a_1,a_2,\underline{a_1+a_2,a_4=a_1+2a_2}
$
\end{itemize}
In the case of (a), we get $\bfa=1,2,3,4$.
In the case of (b), we get $\bfa=(a,2a,b,a+b)$ where $(a,b)=(a_1,a_3)$.
Thus we obtain $\bfa=1,2,3,4$ from (v) of Lemma A.
In the case of (c), we get
$\bfa=(a,b,a+b,2a+b)$ where $(a,b)=(a_1,a_2)$.
Then we obtain $\bfa=1,2,3,4$ from (ii) of Lemma A.
In the case of (d), we get
$\bfa=(a,b,a+b,2b)$ where $(a,b)=(a_1,a_2)$.
Thus we obtain $\bfa=1,2,3,4$ from (iii) of Lemma A.
In the case of (e), we get
$\bfa=(a,b,a+b,a+2b)$ where $(a,b)=(a_1,a_2)$.
Thus we obtain $(a,b)=(2k-1,2k+1)$ for some
$k>0$ from (i) of Lemma A.
Then $S_{\bfa}=S_{2k-1,2k+1,4k,6k+1}$.
Therefore, we conclude that (1) implies (2).
  \end{proof}
  \end{prop}

  \renewcommand{\refname}{References}

     \end{document}